\documentclass[12pt]{amsart}
\usepackage{amsthm,amsfonts, amssymb, amscd}

\newcommand{\C}{{\mathbb C}}

\newcommand{\R}{{\mathbb R}}
\newcommand{\M}{{\mathbb M}}
\newcommand{\Mz}{{\mathbb M}_a}
\newcommand{\T}{{\mathbb T}}

\newcommand{\I}{{\mathbb I}}

\newcommand{\F}{{\mathcal F}}
\newcommand{\G}{{\mathcal G}}
\newcommand{\E}{{\mathcal E}}

\newcommand{\B}{{\mathbb B}}

\newcommand{\FNMR}{{F(N,M;\R)}}
\newcommand{\FNMC}{{F(N,M;\C)}}
\newcommand{\Frames}{{\F[N,M;\R]}}
\newcommand{\FramesC}{{\F[N,M;\C]}}

\newtheorem{thm}{Theorem}
\newtheorem{prop}[thm]{Proposition}

\newtheorem{lemma}[thm]{Lemma}

\newtheorem{cor}[thm]{Corollary}
\theoremstyle{remark}

\def\R{{\mathbb R}}

\newcommand{\OP}[1]{}

\newcommand{\private}[1]{}  

\newcommand{\myprivate}[1]{} 

\newcommand{\norm }[ 1 ]{ \left\| #1 \right\|}
\newcommand{\ip}[2]{ \langle #1, #2 \rangle } 

\newcommand{\ltwoi}{{l^2(\I)}}

\theoremstyle{remark}  

\usepackage[dvips]{epsfig}

\numberwithin{equation}{section}

\begin{document}

\numberwithin{thm}{section}
\title{On signal reconstruction without noisy phase}
\author{Radu~Balan, Pete~Casazza, Dan~Edidin}

\address{\textrm{(R.~Balan)}
Siemens Corporate Research,
755 College Road East,
Princeton, NJ 08540}
\email{radu.balan@siemens.com}

\address{\textrm{(P.~G.~Casazza)}
Department of Mathematics,
University of Missouri,
Columbia, MO 65211}
\email{pete@math.missouri.edu}

\address{\textrm{(Dan~Edidin)}
Department of Mathematics,
University of Missouri,
Columbia, MO 65211}
\email{edidin@math.missouri.edu}
\date{\today}

\thanks{The second author was supported by NSF DMS 0405376
and the third author was supported by NSA MDA 904-03-1-0040}

\begin{abstract}
We construct new classes of Parseval frames for a Hilbert space
which allow signal reconstruction from the absolute value of the
frame coefficients.  As a consequence, signal reconstruction can
be done without using noisy phase or its estimation.  This verifies
a longstanding conjecture of the speech processing community.
\end{abstract}

\maketitle

\section{Introduction}

Reconstruction of a signal using noisy phase or its estimation
can be a critical problem in speech recognition technology.
But, for many years, engineers have believed that speech recognition
should be independent of phase.  By constructing new classes of
Parseval frames for a Hilbert space, we will show that this allows
reconstruction of a signal without using noisy phase or its
estimation.  This verifies the longstanding conjecture of the
speech processing community.

Frames are redundant systems of vectors in a Hilbert spaces.  They
satisfy the
well-known property of perfect reconstruction, in that any vector of the
Hilbert space can be synthesized back from its inner products with the
frame
vectors. More precisely, the linear transformation from the initial
Hilbert
space to the space of coefficients obtained by taking the inner product
of a vector with the
 frame vectors is injective and hence admits a left inverse. This
property
has been succesfully used in a broad spectrum of applications,
including internet coding, multiple antenna coding, optics,
quantum information theory, signal/image processing, and much more.
The purpose of this paper is to study what kind of reconstruction is
possible
if we only have knowledge of the absolute values of the frame
coefficients.

In this paper we consider only finite dimensional frames the reason
being
their direct link to practical applications. Since the same question
can be raised for infinite dimensional frames, we state the problem in
the
setting of abstract frames.

Consider a Hilbert space $H$ with scalar product $\ip{}{}$. A
finite or countable set of vectors $\F=\{f_i;i\in\I\}$ of $H$ is called
a
{\em frame} if there are two positive constants $A,B>0$ such that for
every
vector $x\in H$,
\begin{equation}
\label{eq1.1}
A\norm{x}^2\leq \sum_{i\in I}|\ip{x}{f_i}|^2 \leq B\norm{x}^2
\end{equation}
The frame is {\em tight} when the constants can be chosen equal to one
another, $A=B$. For $A=B=1$, $\F$ is called a {\em Parceval frame}. The
numbers $\ip{x}{f_i}$ are called {\em frame coefficients}.

To a frame $\F$ we associate the {\em analysis} and {\em synthesis
operators}
defined by:
\begin{eqnarray}
T:H\rightarrow \ltwoi & , & T(x)=\{ \ip{x}{f_i} {\}}_{i\in\I} \\
T^*:\ltwoi \rightarrow H & , & T^*(c) = \sum_{i\in I} c_i f_i
\end{eqnarray}
which are well defined due to (\ref{eq1.1}), and are
adjoint to one another. The range of $T$ in $\ltwoi$ is called the {\em
range
of coefficients}. The {\em frame operator} defined by
$S=T^*T:H\rightarrow H$ is invertible by (\ref{eq1.1}) and provides the
perfect reconstruction formula:
\begin{equation}
\label{eq1.2}
x = \sum_{i\in\I}\ip{x}{f_i} S^{-1}f_i
\end{equation}
For more information on frames we refer the reader to \cite{C}.

Consider now the nonlinear mapping
\begin{equation}
\label{eq1.3}
\Mz:H\rightarrow\ltwoi ~,~ \Mz(x) = \{ |\ip{x}{f_i}| {\}}_{i\in\I}
\end{equation}
obtained by taking the absolute value entrywise of the analysis
operator.
Let us denote by $H_r$ the quotient space  $H_r=H/\sim$ obtained by
identifying two vectors that differ by a constant phase factor: $x\sim
y$
 if there is a scalar $c$ with $|c|=1$ so that $y=cx$. For real Hilbert
spaces
$c$ can only be $+1$ or $-1$, and thus $H_r = H/\{\pm 1\}$.
For complex Hilbert spaces $c$ can be any
complex number of modulus one, $c=e^{i\varphi}$, and then $H_r=H/\T^1$,
where
$\T^1$ is the complex unit circle.  In quantum mechanics these
projective rays define quantum states (see \cite{qb}).
Clearly two vectors of $H$ in the same ray would have the same image
through
$\Mz$. Thus the nonlinear mapping $\Mz$ extends to $H_r$ as
\begin{equation}
\M:H_r\rightarrow\ltwoi~,~\M(\hat{x})=\{|\ip{x}{f_i}|{\}}_{i\in\I}~,~
x\in\hat{x}
\end{equation}
The problem we study in this paper is the injectivity of the map $\M$.
When
it is injective, $\M$ admits a left inverse, meaning that any vector
(signal)
 in
$H$ can be reconstructed up to a constant phase factor from the modulus
of
its frame coefficients.

The motivation for this problem comes from two applications in signal
processing, one concerning noise reduction, and the other regarding
 speech recognition. There is also a connection with a problem in
optics that we describe later.

The traditional method of signal enhancement consists of three steps:
first,
the input signal is linearly transformed from its input domain (e.g.
time, or
space) into a transformed domain (e.g. time-frequency, time-scale,
space-scale
 etc.); second, a (nonlinear) estimation operator is applied in this
representation domain; third, a (left) inverse of the linear
transformation
at step one is applied on the signal obtained at step two in order to
synthesize the estimated signal in the input domain. Several linear
transformations have been proposed in the literature and are used in
practice,
 e.g. windowed Fourier transform, wavelet filterbanks, local cosine
basis etc.
 (see \cite{proakis2,vaseghi}). Likewise, many signal estimators have
been
proposed and studied in the literature, some of them statistically
motivated,
 e.g. Wiener (MMSE) filter, Maximum A Posteriori (MAP), Maximum
Likelihood (ML)
etc., others having a rather ad-hoc motivation, e.g. spectral
subtraction,
 psychoacoustically motivated audio and video estimators etc. For
more details see \cite{vanTrees4,proakis,WienerFilter} and many other
books on this topic.  By way of an example let us consider the
Ephraim-Malah
noise reduction method (\cite{EphraimMalah84}) of speech signals. Let
$\{ x(t)~,~t=1,2,\ldots,T\}$ be the samples of a speech signal.
These samples are first transformed into the time-frequency domain
through
\begin{equation}
\label{eq1.7}
X(k,\omega) = \sum_{t=0}^{M-1} g(t)x(t+kN)e^{-2\pi i \omega
\frac{t}{M}}~~,~~
k=0,1,\ldots,\frac{T-M}{N}~,~\omega\in\{0,1,\ldots,M-1\}
\end{equation}
where $g$ is the analysis window, and $M,N$ are respectively the window
size,
and the time step. Next a complicated nonlinear transformation is
applied to
$|X(k,\omega)|$ to produce the MMSE estimate of the short-time spectral
amplitude
\begin{equation}
\label{eq1.x}
Y(k,\omega) =
\frac{\sqrt{\pi}}{2}\frac{\sqrt{v(k,\omega)}}{\gamma(k,\omega)}
exp(-\frac{v(k,\omega)}{2})
[(1+v(k,\omega))I_0(\frac{v(k,\omega)}{2})+v(k,\omega)
I_1(\frac{v(k,\omega)}{2})]|X(k,\omega)|
\end{equation}
where $I_0,I_1$ are modified Bessel functions of zero and first order,
and
$v(k,\omega),\gamma(k,\omega)$ are estimates of certain signal-to-noise
ratios.
The speech signal windowed Fourier coefficients are estimated simply by:

\begin{equation}
\label{eq1.xx}
\hat{X}(k,\omega) = Y(k,\omega) \frac{X(k,\omega)}{|X(k,\omega)|}
\end{equation}
and then are transformed back into time domain through an overlap-add
procedure
\begin{equation}
\label{eq1.10}
\hat{x}(t) = \sum_{k} \sum_{\omega=0}^{M-1} \hat{X}(k,\omega)e^{2\pi i
\omega\frac{t-kN}{M}}h(t-kN)
\end{equation}
where $h$ is the synthesis window. This example illustrates a feature
that is
common to most other signal enhancement algorithms: the nonlinear
estimation in
the representation domain modifies only the amplitude of the transformed
signal, and keeps its noisy phase.  In some applications, such as speech
recognition, reconstruction with noisy phase is a critical problem.
The optimal solution to this problem would occur if we do not need
 the phase at all to perform reconstruction into the input domain.
domain
This paper addresses exactly this issue.

Consider now the problem of automatic speech recognition (ASR) systems.
Given
a voice signal $\{x(t),t=1,2,\ldots,T\}$, the ASR outputs a sequence of
recognized phonemes from an alphabet. Most ASR systems use different
kind of
{\em cepstral } coefficient statistics (see
\cite{rabinerjuang,becchetti})
as described next. The voice signal is transformed into the
time-frequency
domain by the same discrete windowed Fourier transform (\ref{eq1.7}).
The
(real) cepstral coefficients $C_x(k,\omega)$ are defined as the
logarithm of
the modulus of $X(k,\omega)$:
\begin{equation}
\label{eq1.11}
C_x(k,\omega)={\rm log}(|X(k,\omega)|)
\end{equation}
There are two rationals for using this object. First note the recorded
signal
$x(t)$ is a convolution of the voice signal $s(t)$ with the
source-to-microphone (channel) impulse response $h$. In the
time-frequency domain,
convolution
becomes (almost) multiplication, and the cepstral coefficients decouple
\begin{equation}
\label{eq1.12}
C_x(k,\omega)={\rm log}(|H(\omega)|) + C_s(k,\omega)
\end{equation}
where $H(\omega)$ is the channel transfer function, and $C_s$ is the
voice
signal cepstral coefficient. Since the channel transfer function is
time-invariant, by subtracting the time average we obtain
\begin{equation}
\label{eq1.13}
F_x(k,\omega) = C_x(k,\omega) - \E[C_x(\cdot,\omega)] = C_s(k,\omega) -
\E[C_s
(\cdot,\omega)]
\end{equation}
where $\E$ is the time average operator. Thus $F_x$ encodes information
about
the speech signal alone, independent of the reverberant environment.

The second reason for using $C_x$, and thus $F_x$, is the widespread
belief
in the speech processing community that phase does not matter in speech
recognition. Hence, by taking the modulus in (\ref{eq1.11}) one does not
lose information about the message (nor the messanger, as in some
speaker
identification algorithms).

Returning to the ASR system, the corrected cepstral coefficients $F_x$
are
fed into several hidden Markov models (HMMs), one HMM for each phoneme.
The
outputs of these HMMs give the utterance likelihood of a particular
phoneme.
 The ASR system then chooses the phoneme with the largest likelihood.

In the two classes of signal processing algorithms described above the
transformed domain signal either has a secondary role,
or has none whatsoever. This observation led us to consider the
information
loss introduced by taking the modulus of a redundant representation.
Clearly a constant phase is always lost, however is this the only loss
of information with respect to the original signal? This is the problem
we analyze in this paper.

There is also a closely connected problem in optics with applications to
X-ray, crystallography, electron microscopy, and coherence theory
 see \cite{fienup78,fienup82,reviewBates,liu90}.
This problem is to
reconstruct a discrete signal from the modulus of its Fourier transform
under constraints in both the original and the Fourier domain. For
finite
signals the approach uses the Fourier transform with redundancy 2. All
signals with the same modulus of the Fourier transform satisfy a
polynomial
factorization equation.
In dimension one this factorization has an exponential number of
possible solutions. In higher dimensions the factorization is shown to
have
generically a unique solution (see \cite{hayes82}).

The organization of the paper is as follows. Section \ref{sec2}
presents the analysis of real frames;  section \ref{sec3}
analyzes the case of complex frames. 

\section{Analysis of $\M$ for Real Frames\label{sec2}}

Consider the case $H=\R^N$, and the index set $\I$ has cardinality $M$,
 $\I=\{1,2,\ldots,M\}$. Then $\ltwoi \simeq \R^M$.

For a frame $\F=\{f_1,\ldots,f_M\}$ of $\R^N$ we denote by $T$ the
analysis
operator,
\begin{equation}
\label{eq2.1}
T:\R^N\rightarrow\R^M ~~,~~ T(x)=\sum_{k=1}^M
\ip{x}{f_k}e_k~~,~~x\in\hat{x}
\end{equation}
where $\{e_1,\ldots,e_M\}$ is the canonical basis of $\R^M$. We let $W$
denote
the range of the analysis map $T(\R^N)$ which is an $N$-dimensional
subspace of
$\R^M$.  Recall the nonlinear map we are interested in is
\begin{equation}
\label{eq2.2}
\M^{\F}:\R^N/\{\pm 1\}\rightarrow\R^M~~,~~
\M^{\F}(\hat{x})=\sum_{k=1}^M |\ip{x}{f_k}|e_k
\end{equation}
When there is no danger of confusion, we shall drop $\F$ from the
notation.

Two frames $\{f_i \}_{i\in I}$ and $\{g_i \}_{i\in I}$ 
are {\it equivalent} if there is an
invertable
operator $T$ on $H$ with $T(f_i ) = g_i$, for all $i\in I$.
It is known that two frames are equivalent
if and only if their associated analysis operators have the same
range (see \cite{radu98c,hanlars00}). We deduce that $M$-element frames
on
$\R^N$ are parametrized by the fiber bundle $\FNMR$, which is
the $GL(N,\R)$-bundle over the Grassmanian $Gr(N,M)$.

First we reduce our analysis to equivalent classes of frames:

\begin{prop}\label{P2.1}
For any two frames $\F$ and $\G$ that have the same range of
coefficients,
$\M^{\F}$ is injective if and only if $\M^{\G}$ is injective.
\end{prop}

\begin{proof}
Any two frames $\F=\{f_k\}$ and $\G=\{g_k\}$ that have the same range of
coefficients are equivalent, i.e. there is an invertible
$R:\R^N\rightarrow\R^N$ so that $g_k=Rf_k$, $1\leq k\leq M$.
Their associated nonlinear maps $\M^{\F}$, and respectively $\M^{\G}$,
satisfy $\M^{\G}(x)=\M^{\F}(R^*x)$. This shows that $\M^{\F}$ is
injective if and only if $\M^{\G}$ is injective. Consequently the
property of
injectivity of $\M$ depends only on the subspace of coefficients $W$
in $Gr(N,M)$. 
\end{proof}

This result says that for two frames corresponding to two points in the
same fiber of $\FNMR$, the injectivity of their associated nonlinear
maps would jointly hold true or fail. Because of this result we shall
always assume the induced topology by the base manifold $Gr(N,M)$ of the
fiber
bundle $\FNMR$ into the set of $M$-element frames of $\R^N$.

If $\{f_i \}_{i\in I}$ is a frame with frame operator $S$ then
$\{S^{-1/2}f_i\}_{i\in I}$ is a Parseval frame which is equivalent to
$\{f_i \}_{i\in I}$ and called the {\it canonical Parseval frame}
associated to $\{f_i \}_{i\in I}$.  Also, $\{S^{-1}f_i\}_{i\in I}$ is
a frame equivalent to $\{f_i \}_{i\in I}$ and is called the
{\it canonical dual frame} associated to $\{f_i \}_{i\in I}$.
Proposition 2.1 shows that when the nonlinear map $\M^{\F}$ is injective
then the same property holds for the canonical dual frame and
the canonical Parseval frame.

Given $S \subset \{1, \ldots , M\}$ define a map
$\sigma_S\colon \R^M \to \R^M$ by the formula
$$
\sigma_S(a_1, \ldots , a_M)
= ((-1)^{S(1)}a_1, \ldots , (-1)^{S(M)}a_M).
$$
 Clearly
$\sigma^2_S = id$ and $\sigma_{S^\complement} = -\sigma_S$
 where $S^\complement$ is the complement of $S$. We let $L^S$ denote the
$|S|$-dimensional linear subspace of $\R^M$ where
$L^{S} = \{(a_1, \ldots , a_M) | a_i = 0, i \in S\}$, and we let
$P_S:\R^M\rightarrow L^S$ denote the orthogonal projection onto this
subspace.
Thus $(P_S(u))_i=0$, if $i\in S$, and $(P_S(u))_i=u_i$, if
$i\in S^\complement$. For every vector $u\in\R^M$, $\sigma_S(u)=u$ iff
$u\in L^S$. Likewise $\sigma_S(u)=-u$ iff $u\in L^{S^\complement}$. Note

\[ P_S(u) =
\frac{1}{2}(u+\sigma_S(u))~~~,~~~P_{S^\complement}(u)=\frac{1}{2}
(u-\sigma_S(u)) \]

\begin{thm}[Real Frames]\label{T2.1}
If $M \geq 2N-1$ then for a generic frame $\F$, $\M$ is injective.
\end{thm}
By {\it generic} we mean an open dense subset of the set of all
$M$-element
frames in $\R^N$.
\begin{proof}
Suppose that $x$ and $x'$ have the same image under $\M=\M^{\F}$.
Let $a_1, \ldots a_M$ be the frame coefficients
of $x$  and $a'_1, \ldots a'_M$ the frame coefficients for $x'$.
Then $a'_i = \pm a_i$ for each $i$. In particular there is a subset $S
\subset
\{1, \ldots , M\}$ of indices such that $a'_i = (-1)^{S(i)}a_i$
where the function
$S(i)$ is the characteristic function of $S$ and is defined by the rule
that $S(i) = 1$ if $i \in S$ and $S(i) = 0$
if $i \notin S$. Then two vectors $x$, $x'$ have the same
image under $\M$ if and only there is a subset
$S \subset \{1, \ldots , M\}$ such that
$(a_1, \ldots a_M)$ and $((-1)^{S(1)}a_1, \ldots , (-1)^{S(M)}a_M)$
are both in $W$ the range of coefficients associated to $\F$.

To finish the the proof we will show that when $M \geq 2N-1$ such a
condition is impossible for a generic subspace $W \subset \R^N$.
This means that the set
of such $W$'s is a dense (Zariski) open set in the Grassmanian
$Gr(n,m)$. In particular the probability that a randomly chosen $W$ will
satisfy this condition is $0$.

To finish the proof of the theorem we need the following lemma.
\begin{lemma}\label{L2.1}
If $M \geq 2N-1$ then the following holds for a generic $N$-dimensional
subspace $W\subset \R^M$. Given $u\in W$ then $\sigma_S(u) \in W$
iff $\sigma_S(u) = \pm u$.
\end{lemma}

\begin{proof}[Proof of the Lemma]
Suppose $u \in W$ and $\sigma_S(u) \neq \pm u$ but
$\sigma_S(u) \in W$. Since $\sigma_S$ is an involution,
$u + \sigma_S(u)$ is fixed by $\sigma_S$ and is non-zero.
Thus $W \cap L^{S} \neq 0$. Likewise
$$0\neq u - \sigma_S(u) =
u + \sigma_{S^\complement}(u).$$ Hence
$W \cap L^{S^\complement} \neq 0$.

Now $L^S$ and $L^{S^c}$ are fixed linear subspaces of dimension $M -
|S|$
and $|S|$. If $M \geq  2N-1$ then one of these subspaces has
codimension greater than or equal to $N$.
However a generic linear subspace $W$ of dimension
$N$ has $0$ intersection with a fixed linear subspace of codimension
greater
than or equal to $N$.
Therefore, if $W$ is generic and $x,\sigma_S(x) \in W$ then
$\sigma_S(x) = \pm x$ which ends the proof of Lemma.
\end{proof}

The proof of the theorem now follows from the fact that if
$W$ is in the intersection of generic conditions imposed
by the proposition for each subset $S \subset \{1, \ldots , M\}$
then $W$ satisfies the conclusion of the theorem.
\end{proof}

Note what the above proof actually shows:

\begin{cor}\label{C2}
The map $\M$ is injective if and only if whenever there is a non-zero
element $u\in W \subset \R^M$ with $u\in L^S$, then $W\cap
L^{S^\complement}
= \{0\}$.
\end{cor}

Next we observe that this result is best possible.

\begin{prop}\label{P1}
If $M\le 2N-2$, then the result fails for all $M$-element frames.
\end{prop}

\begin{proof}
Since $M\le 2N-2$, we have that $2M-2N+2 \le M$.
Let $(e_i )_{i=1}^{M}$ be the canonical
orthonormal basis of $\R^M$.  We can write $(e_i)_{i=1}^{M} =
(e_i)_{i=1}^{k}
\cup (e_i)_{i=k+1}^{M}$ where both $k$ and $M-k$ are $\ge M-N+1$.

Let $W$ be any $N$-dimensional subspace of $\R^M$.
Since dim $W^{\perp}=M-N$, there exists a nonzero vector
$u\in$span $\{e_{i}\}_{i=1}^{k}$ so that $u\perp W^{\perp}$, hence $u\in
W$.
Similarly, there is a nonzero vector $v$ in span $\{e_{i}\}_{i=k+1}^{M}$
with $v\perp W^{\perp}$, that is $v\in W$.
By the above corollary, $\M$ cannot be injective. In fact
$\M(u+v)=\M(u-v)$.
\end{proof}

The next result gives an easy way for frames to satisfy the condition
above.

\begin{cor}\label{C1}
If $\F$ is a $M$-element frame for $\R^N$ with
$M\ge 2N-1$ having the property that every $N$-element subset of the
frame is linearly independent, then $\M$ is injective.
\end{cor}

\begin{proof}
Given the conditions, it follows that $W$ has no elements which are
zero in $N$ coordinates and so the Corollary holds.
\end{proof}

\begin{cor}\label{C2.3}\mbox{}

\begin{enumerate}
\item If $M=2N-1$, then the condition given in Corollary \ref{C1} is
also
necessary.
\item If $M\ge 2N$, this condition is no longer necessary.
\end{enumerate}
\end{cor}

\begin{proof}
(1) For the first part we will prove the contrapositive.
Let $M=2N-1$ and assume there is an $N$-element subset $(f_i)_{i\in S}$
of $\F$ which is not linearly independent. Then there is a non-zero
$x\in(span(f_i)_{i\in S})^{\perp}\subset\R^N$.  Hence,  $0\neq u=T(x)\in
L^S\cap W$. On the other hand, since $dim(span(f_i)_{i\in
S^\complement})\leq
N-1$, there is a non-zero $y\in(span(f_i)_{i\in S^\complement})^{\perp}
\subset\R^N$ so that $0\neq v=T(y)\in L^{S^\complement}\cap W$.
Now, by Corollary \ref{C2}, $\M$ is not injective.

(2) If $M\ge 2N$ we construct an $M$-element frame for $\R^N$ that has
an
$N$-element linearly dependent subset.  Let
$\F'=\{f_1,\ldots,f_{2N-1}\}$ be
a frame for $\R^N$ so that any $N$-element subset is linearly
independent.
By Corollary \ref{C2}, the map $\M^{\F'}$ is injective. Now extend this
frame
to $\F=\{f_1,\ldots,f_M\}$ by $f_{2N}=\cdots=f_M=f_{2N-1}$. The map
$\M^\F$
extends $\M^{\F'}$ and therefore remains injective, whereas clearly any
$N$-element subset that contains two vectors from
$\{f_{2N-1},f_{2N},\ldots,
f_{M}\}$ is no longer linearly independent.
\end{proof}

{\it Remark}: The frames above can easily be constructed ``by hand''.
Start with an orthonormal basis for $\R^N$, say $(f_i)_{i=1}^{N}$.
Assume we have constructed sets of vectors $(f_i)_{i=1}^{M}$ such that
every subset of N vectors is linearly independent.  Look at the span
of all of the (N-1)-element subsets of $(f_i)_{i=1}^{M}$.  Pick
$f_{M+1}$ not in the span of any of these subsets.  Then
$(f_i)_{i=1}^{M+1}$ has the property that every N-element subset is
linearly independent.

Now we will give a slightly different proof of this
result which gives necessary and sufficient conditions
for a frame to have the required properties.

\begin{thm}\label{T2.2}
Let $(f_i)_{i=1}^{M}$ be a frame for ${\R^N}$.
The following are equivalent:

(1)  The map $\M$ is injective.

(2)  For every subset $S\subset \{1,2,\ldots , M\}$,
either $\{f_{i}\}_{i\in S}$ spans $\R^N$ or
$\{f_{i}\}_{i\in S^c}$ spans $\R^N$.
\end{thm}

\begin{proof}
$(1)\Rightarrow (2)$:  We prove the contrapositive.
So assume that there is a subset $S\subset
\{1,2,\ldots M \}$ so that neither
$\{f_i~;~i\in S\}$ nor $\{f_i~;~i\in S^\complement\}$ spans $\R^N$.
Hence there are non-zero vectors $x,y\in\R^N$ so that $x\perp
span(f_i)_{i\in S}$ and $y\perp span(f_i)_{i\in S^\complement}$. Then
$0\neq T(x)\in L^S\cap W$ and $0\neq T(y)\in L^{S^\complement}\cap W$.
Now
by Corollary \ref{C2} we have that $\M$ cannot be injective.

$(2) \Rightarrow (1)$:  Suppose $\M(\hat{x})=\M(\hat{y})$ for some
$\hat{x},\hat{y}\in\R^N/\{\pm 1\}$. This means for every $1\le j\le M$,
$|\ip{x}{f_j}|=|\ip{y}{f_j}|$ where $x\in\hat{x}$ and $y\in\hat{y}$.
Let
\begin{equation}
S = \{j~:~\ip{x}{f_j} = -\ip{y}{f_j} \}.
\end{equation}
Note
\begin{equation}
S^\complement = \{j~:~\ip{x}{f_j} = \ip{y}{f_j} \}
\end{equation}
Now, $x+y \perp span(f_i)_{i\in S}$ and $x-y\perp span(f_i)_{i\in
S^\complement}$. Assume that $\{f_i~;~i\in S\}$ spans $\R^N$. Then
$x+y=0$
and thus $\hat{x}=\hat{y}$. If $\{f_i~;~i\in S^\complement\}$ spans
$\R^N$
then $x-y=0$ and again $\hat{x}=\hat{y}$. Either way $\hat{x}=\hat{y}$
which
proves $\M$ is injective. 
\end{proof}

For $M<2N-1$ there are plenty of frames for which $\M$ is not injective.
However for a generic frame, we can show the set of rays that can be 
reconstructed from the image under $\M$ is open dense in $\R^N/\{\pm 1\}$.

\begin{thm}
Assume $M>N$. Then for a generic frame $\F\in\Frames$, the set of vectors 
$x\in\R^N$ so that $(\M^\F)^{-1}(\M_a^\F(x))$ consists of one point in 
$\R^N/\{\pm 1\}$ has dense interior in $\R^N$.
\end{thm}

\begin{proof}

Let $\F$ be a $M$-element frame in $\R^N$. Then $\F$ is similar to
a frame $\G$ which consists of the union of the canonical basis of $\R^N$, 
$\{d_1,\ldots,d_N\}$, with some other set of $M-N$ vectors. 
Let $\G=\{g_k~;~1\leq k\leq M\}$. Thus $g_{k_j}=d_j$, $1\leq j\leq N$, for 
some $N$ elements $\{k_1,k_2,\ldots,k_N\}$ of $\{1,2,\ldots,M\}$.
Consider now the set $\B$ of frames $\F$ so that its similar frame $\G$ 
constructed above has a vector $g_k$ with all entries non-zero,
\[ \B = \{ \F\in\Frames~|~\F\sim\G=\{g_k\}~,~\{d_1,\ldots,d_N\}\subset\G~,~
\prod_{j=1}^N\ip{g_{k_0}}{d_j}\neq 0,{\rm for~some}~k_0\} \]
Clearly $\B$ is open dense in $\Frames$. Thus  generically $\F\in\B$. Let
$\G$ be its similar frame satisfying the condition above.
We want to prove the set $X=X^\F$ of vectors $x\in\R^N$ so that 
$(\M^\G)^{-1}(\M_a^\G(x))$ has more than one point is {\em thin}, i.e. it is 
included into a set whose complement is open and dense in $\R^N$. 
We claim $X\subset \cup_{S}(V_S^+\cup V_S^-)$ where 
$(V_S^{\pm})_{S\subset\{1,2,\ldots,N\}}$ are linear 
subspaces of $\R^N$ of codimension $1$ indexed by subsets $S$ of 
$\{1,2\ldots,N\}$.
This this claim will conclude the proof of Theorem.

To verify the claim,
let $x,y\in\R^N$ be so that $\M_a^\G(x)=\M_a^\G(y)$ and yet $x\neq y$, nor 
$x\neq -y$. Since $\G$ contains the canonical basis of $\R^N$, $|x_k|=|y_k|$
for all $1\leq k\leq N$. Then there is a subset $S\subset\{1,2,\ldots,N\}$
so that $y_k=(-1)^{S(k)}x_k$. Note $S\neq\emptyset$, nor 
$S\neq\{1,2,\ldots,N\}$. Denote by $D_S$ the diagonal $N\times N$ matrix
$(D_S)_{kk}=(-1)^{S(k)}$. Thus $y=D_Sx$, and yet $D_S\neq \pm I$. Let
$g_{k_0}\in\G$ be so that none of its entries vanishes. Then 
$|\ip{x}{g_{k_0}}|=|\ip{y}{g_{k_0}}|$ implies 
\[ \ip{x}{(I\pm D_S)g_{k_0}} =0 \]
This proves the set $X^\G$ is included into the union of $2(2^N-2)$
linear subspaces of codimension $1$,
\[ \cup_{S\neq\emptyset,S^{\complement}\neq\emptyset}\{(I-D_S)g_{k_0}{\}}^{\perp}
\cup \{(I+D_S)g_{k_0} {\}}^{\perp} \]
Since $\F$ is similar to $\G$, $X^\F$ is included into the image of the 
above set through a linear invertible map, which proves the claim. 

\end{proof}

\section{Analysis of $\M$ for Complex Frames\label{sec3}}

In this section the Hilbert space is $\C^N$. For an $M$-element frame
$\F=\{f_1,\ldots,f_M\}$ of $\C^N$ the analysis operator is defined by
(\ref{eq2.1}), where the scalar product is $\ip{x}{y}=\sum_{k=1}^N
x_k\overline{y_k}$. The range of coefficients, i.e. the range of the
analysis operator, is a complex $N$-dimensional
subspace of $\C^M$ that we denote again by $W$. The nonlinear map we are
studying is given by
\begin{equation}
\label{eq3.1}
\M^\F:\C^N/\T^1\rightarrow \C^M~~,~~\M^\F(\hat{x})=\sum_{k=1}^M
|\ip{x}{f_k}|e_k~~,~~x\in\hat{x}
\end{equation}
where two vectors $x,y\in\hat{x}$ if there is a scalar $c\in\C$ with
$|c|=1$ so that $y=cx$.

By the equivalence results proved in \cite{radu98c,hanlars00} we obtain
that $M$-frames of $\C^N$ are parametrized by points of the fiber bundle
$\FNMC$, the $GL(N,\C)$-bundle over the complex Grassmanian
$Gr(N,M)^\C$.

Proposition \ref{P2.1} holds true for complex frames as well.
Thus without loss of generality we shall work with the topology induced
by the
base manifold of $\FNMC$ into the set of $M$-element frames of $\C^N$.

As in the real case we reduce the question about $M$-element frames in
$\C^N$ to a question about the Grassmanian of $N$-planes in $\C^M$.
First we
prove the following

\begin{thm}\label{T3.1}
If $M \geq 4N -2$ then the generic $N$-plane
$W$ in $\C^M$ has the property
that if $v = (v_1, \ldots,  v_M)$ and $w = (w_1, \ldots ,w_M)$
are vectors in $W$ such that $|v_i| = |w_i|$ for all $i$ then
$v = \lambda w$ for some complex number $\lambda$ of modulus $1$.
\end{thm}

\begin{proof}
We will say that an $N$-plane $W$ has {\em property (*)} if
there are non-parallel vectors $v,w$ in $W$ such that $|v_i| = |w_i|$
for all $i$. Recall two vectors $x,y$ are {\em parallel} if there
is a scalar $c\in\C$ so that $y=cx$.

Given a $N$-plane $W$ we may assume,
after reordering the coordinates on $\C^M$, that
$W$ is the span of the rows of a $N\times M$ matrix of the form

$$\left[ \begin{array}{ccccccc}
1 & 0 & \ldots & 0 & u_{N+1,1} & \ldots & u_{M,1}\\
0 & 1 & \ldots & 0 & u_{N+1,2} & \ldots & u_{M,2}\\
| & | &    |   & | & |         &  |     &        \\
0 & 0 & \ldots & 1 & u_{N+1,N} & \ldots & u_{M,N}
\end{array}
\right]$$
where the $N(M-N)$ entries $\{u_{i,j}\}$ are viewed as indeterminates.
Thus $Gr(N,M)^\C$ is isomorphic to $\C^{N(M-N)}$ in a neighborhood of
$W$.

Now suppose that $W$ satisfies (*) and $v$ and $w$ are two non-parallel
vectors
whose entries have the same modulus. Our choice of basis for
$W$ ensures that one of the first $N$ entries in $v$ (and hence $w$)
are non-zero. Since we only care about these vectors up to rescaling
we may assume, after reordering, that $v_1 = w_1 = 1$. Also the
vectors are assumed non-parallel so we may assume that
$v_i \neq w_i \neq 0$ for some $i \leq N$. After yet again
reordering we can assume that $v_2 \neq w_2  \neq 0$.

Set $\lambda_1 = 1$. By assumption there are numbers
$\lambda_2, \ldots , \lambda_M \in \T^1$ with $\lambda_2\neq 1$
such that $w_i = \lambda_i v_i$ for $i = 1, \ldots,  M$. Expanding
in terms of the basis for $W$ we have
for $i > N$,  $v_i = \sum_{j =1}^N v_j u_{i,j}$ and
$w_i = \sum_{j =1}^N \lambda_j v_j u_{i,j}$.
Thus if $W$ satisfies (*) there must be $\lambda_2, \ldots, \lambda_N
\in \T^1$ (with $\lambda_2 \neq 1$) and $v_2, \ldots v_N \in \C$ such
that
for all $N+1\le i \le M$ we have

\begin{equation} \label{eq.goodeqns}
|\sum_{j =1}^N v_j u_{i,j}| = |\sum_{j =1}^N \lambda_j v_j u_{i,j}|.
\end{equation}

Consider the variety $Y$ of all tuples
$$(W, v_2, \ldots , v_N, \lambda_2, \ldots , \lambda_N)$$
as above. Since $v_2 \neq 0$ and $\lambda_2 \neq 1$ this variety
is locally isomorphic to the real $2N(M-N)+ 3N -3$-dimensional
variety
$$\C^{N(M-N)} \times (\C \smallsetminus\{0\}) \times (\C)^{N-2} \times
(\T^1\smallsetminus\{1\})\times
(\T^1)^{N-2}.
$$
 The locus in $Gr(N,M)^\C$ of planes satisfying property (*) is denoted
by $X$.
 This variety is the image under projection to the first factor of $Y$
cut out by the $M-N$ equations (\ref{eq.goodeqns}) for $N+1\le i\le M$.
The analysis of these equations is summarized by the following result.

\begin{lemma} \label{l3.2}The $M-N$ equations in \eqref{eq.goodeqns}
are independent. Hence $X$ is a variety of real dimension at most
$2N(M-N)+ 3N -3 - (M-N)$
\end{lemma}
\begin{proof}[Proof of Lemma 3.2]

For any choice of $0\neq v_2, v_3, \ldots v_N$ and $1\neq \lambda_2,
\lambda_3, \ldots \lambda_N$ the equation
$$
|\sum_{j =1}^M v_j u_{i,j}|^2 = |\sum_{j =1}^M \lambda_j v_j
u_{i,j}|^2
$$
is non-degenerate. Since the variables $u_{i,1}, \ldots , u_{i,N}$
appear in exactly one equation, these equations (for fixed
$v_2, v_3, \ldots , v_N, \lambda_2, \ldots , \lambda_N$) define
a subspace
of $\C^{N(M-N)}$ of real codimension at least $M-N$.
Since this is true for all choices, it follows that the equations
are independent.
\end{proof}

From this lemma it follows that the locus of $N$-planes satisfying
(*) has (local) real dimension $2N(M-N) + 3N -3 - (M-N)$. Therefore
if $3N -3 - (M-N) < 0$, i.e. if $M \geq 4N-2$, this locus can not be all
of $Gr(N,M)^\C$. This ends the proof of Theorem \ref{T3.1}.
\end{proof}

The main result in the complex case then follows from Theorem
\ref{T3.1}.

\begin{thm}[Complex frames] \label{thm.complex}
If $M \geq 4N -2$ then $\M^\F$ is injective for a generic frame
 $\F=\{f_1, \ldots f_N\}$.
\end{thm}

Lemma \ref{l3.2} yields the following result.
\begin{thm}\label{t3.4} If $M\geq 2N$ then for a generic frame $\F\in\FramesC$
the set of vectors $x\in\C^N$ such that $(\M^\F)^{-1}(\M_a^\F(x))$ has one
point in $\C^N/\T^1$ has dense interior in $\C^N$.
\end{thm}

\begin{proof}
By Lemma \ref{l3.2}, for a generic frame the $M-N$ equations 
(\ref{eq.goodeqns}) in $2(N-1)$ indeterminates 
$(v_2,\ldots,v_N,\lambda_2,\ldots,\lambda_N)$ are independent.
Note there are $3(N-1)$ real valued unknowns and $M-N$ equations.
Hence the set of $\{(v_2,\ldots,v_N)\}$ in $\C^{N-1}$ for which there are
$(\lambda_2,\ldots,\lambda_N)$ such that (\ref{eq.goodeqns}) has solution
in $(\C \smallsetminus\{0\}) \times (\C)^{N-2} \times
(\T^1\smallsetminus\{1\})\times
(\T^1)^{N-2}$ has real dimension at most $3(N-1)-(M-N)=4N-3-M$. For $M\geq 2N$
 it follows $3(N-1)-(M-N)<2(N-1)$ which shows the set of $v=(v_1,\ldots,v_N)$
 such that $(\M^\F)^{-1}(\M_a^\F(v))$ has more than one point is thin in
$\C^N$, i.e. its complement has dense interior.
\end{proof}

We do not know the precise optimal bound for the complex case but
we believe it is $4N-2$.  However, this case is different from the real
case in that complex frames with only $2N-1$ elements cannot have
$\M^{\F}$ injective.  To see this we observe that the proof of Theorem
2.8
$(1)\Rightarrow (2)$ does not use the fact that the frames are real.  So
in the complex case we have:

\begin{prop}
If $\{f_j \}_{j\in I}$ is a complex frame and $\M^{\F}$ is injetive,
then
for every $S\subset \{1,2, \ldots , M\}$, if $L^S \cap W \not= \{0\}$
then $L^{S^c} \cap W = \{0\}$.  Hence, for every such $S$,
either $\{f_j\}_{j\in S}$ or $\{f_j \}_{j\in S^c}$ spans $H$.
\end{prop}

Now we can show that complex frames must contain at least
$2N$-elements for $\M^{\F}$ to be injective.

\begin{prop}
[Complex frames]
If $\M^{\F}$ is injective then $M\ge 2N$.
\end{prop}

\begin{proof}
We assume that $M=2N-1$ and show that in this case $\M^{\F}$ is
not injective.  Let $\{z_j \}_{j=1}^{N}$ be a basis for $W$ and let
$P$ be the orthogonal projection onto the first $N-1$ unit vectors
in $\C^M$.  Then $\{Pz_j \}_{j=1}^{N}$ sits in an $N-1$-dimensional
space and so there are complex scalars $\{a_j \}_{j=1}^{N-1}$,
not all zero, so that $\sum a_j Pz_j = 0$.  That is, there is a vector
$0\not= y\in W$ with support $y\subset \{N,N+1 , \ldots 2N-1\}$.
Similarly, there is a vector $0\not= x\in W$ with support
$x\subset \{1,2,\ldots , N\}$.  If $x(N) = 0$ or $y(N) = 0$ we
contradict Proposition 3.4.  Also, if $x(i) = 0$ for all $i<N$, then
$(y-cx)(N) = 0$ for $c = y(N)\frac{\overline{x(N)}}{|x(N)|^2}$.
Now, $x,y-cx$ are in $W$ and have disjoint support so our
map is not injective.  Otherwise, let
$$
z = \frac{\overline{x(N)}}{|x(N)|^2},\ \ w = i\frac{\overline{y(N)}}
{|y(N)|^2}.
$$
Now, $z,w\in W$ and $z(N) = 1$ and $w(N) = i$.  Hence,
$|z+w| = |z-w|$.  It follows that there is a complex number
$|c|=1$ so that $z+w = c(z-w)$.  Since $z_i \not= 0$ for some
$i<N$ we have that $c=1$ and $w=0$ which is a contradiction.
\end{proof}

\bibliographystyle{plain}


\end{document}